\documentclass{amsart}

\usepackage{graphicx,tikz,amssymb}
\usetikzlibrary{arrows}
\usepackage[justification=centering]{caption}

\newtheorem{theorem}{Theorem}[section]

\newtheorem{corr}[theorem]{Corollary}

\DeclareMathOperator{\olws}{olws}
\DeclareMathOperator{\olwi}{olwi}

\author{Csaba Bir\'o}
\email{csaba.biro@louisville.edu}

\author{Israel R. Curbelo}
\email{israel.curbelo@louisville.edu}

\address{Department of Mathematics, University of Louisville, Louisville, KY 40292}

\title[On-Line Chain Partitioning of Semi-Orders]{Improved lower bound on the on-line chain partitioning of semi-orders with representation}

\begin{document}

\begin{abstract}
An on-line chain partitioning algorithm receives a poset, one element at a time, and irrevocably assigns the element to one of the chains in the partition. The on-line chain partitioning problem involves finding the minimal number of chains needed by an optimal on-line algorithm. Chrobak and \'Slusarek considered variants of the on-line chain partitioning problem in which the elements are presented as intervals and intersecting intervals are incomparable. They constructed an on-line algorithm which uses at most $3w-2$ chains, where $w$ is the width of the interval order, and showed that this algorithm is optimal. They also considered the problem restricted to intervals of unit-length and while they showed that first-fit needs at most $2w-1$ chains, over $30$ years later, it remains unknown whether a more optimal algorithm exists. In this paper, we improve upon previously known bounds and show that any on-line algorithm can be forced to use $\lceil\frac{3}{2}w\rceil$ chains to partition a semi-order presented in the form of its unit-interval representation. As a consequence, we completely solve the problem for $w=3$.
\end{abstract}

\maketitle

\section{Introduction}

An on-line chain partitioning algorithm receives a poset $(X,P)$ in the order of its elements $x_1,\ldots,x_n$ and constructs an on-line chain partition. This means that the chain to which the element $x_i$ is assigned to depends solely on the subposet induced by the elements $\{x_1,\ldots,x_{i-1}\}$ and on the chains to which they were assigned to. The effiency of an algorithm is measured with respect to the minimum number of chains needed by an optimal off-line algorithm. By Dilworth's theorem, a poset of width $w$ can always be partitioned off-line into $w$ chains. However, this is not the case when the poset is presented in an on-line manner.

The on-line width $\text{olw}(w)$ of the class of posets of width at most $w$ is the largest integer $k$ for which there exists a strategy that forces any algorithm to use $k$ chains to partition a poset of width $w$. The exact value of $\text{olw}(w)$ remains unknown for $w>2$. Kierstead \cite{kie-81} was the first to prove that $\text{olw}(w)$ was bounded. He constructed an on-line algorithm which uses at most $(5^w-1)/4$ chains to partition a poset of width $w$. Nearly 30 years later, Bosek and Krawczyk \cite{bos-kra-15} presented the first on-line algorithm that uses subexponentially many chains. An easier and more efficient on-line algorithm was presented in \cite{bos-kie-18} by Bosek, Kierstead, Krawczyk, Matecki, and Smith. Most recently, Bosek and Krawcyk \cite{bos-kra-21} constructed an on-line algorithm which needs at most $w^{O(\log \log w)}$ chains to partition a poset of width $w$.  On the other hand, Szemer\'edi provided an unpublished argument (see \cite{kie-86} for a proof) which shows that any algorithm could be forced to use $\binom{w+1}{2}$ chains to partition a poset of width $w$. Szemer\'edi's arguement was later improved in a survey paper by Bosek, Felsner, Kloch, Krawczyk, Matecki and Micek \cite{survey} where they improved the lower bound to $(2-o(1))\binom{w+1}{2}$. Many variants have branched from the general problem by restricting the class of posets further or by restricting the way the poset is presented. We refer the reader to the survey paper \cite{survey} for an overview of problems and results in this field.

%The problem restricted to the class of interval orders is equivalent to the on-line coloring of interval graphs with maximum clique-size $w$ and was solved by Kierstead and Trotter \cite{kie-tro-81} in the language of recursive combinatorics. Chrobak and \'Slusarek \cite{chr-slu-88} considered the problem in regard to dynamic storage allocation. In this setting, elements are presented in the form of intervals instead of points. They independently arrived at the same optimal efficiency. In the same paper, Chrobak and \'Slusarek consider the problem restricted to intervals of unit-length and determine the problem to be harder to analyze than the general problem. They show that the efficiency of First Fit Algorithm is $2w-1$ providing 

In this paper, we focus on the on-line width $\text{olws}_R(w)$ of the class of semi-orders presented in the form of a unit-interval representation and improve upon previously known bounds by proving the following result.
\begin{theorem}
    $\olws_R(w)\geq \lceil\frac{3}{2}w\rceil$ for $w>1$.
\end{theorem}

In the following section we provide some background and a brief sketch of previous results for closely related problems. While some of the problems have been considered in a graph-theoretic setting, as in \cite{gya-leh-88} by Gy\'arf\'as and Lehel, we translate all results to the poset-theoretic setting in order to emphasize the connections between them and for simplicity. In section $3$, we provide a proof of the main theorem.

\section{Posets and On-Line Width}

A poset $(X,P)$ is an interval order if there is a function $I$ which assigns to each element $x\in X$ a closed interval $I(x)=[l_x,r_x]$ on the real line so that for all $x_1,x_2\in X$ we have $x_1<x_2$ if and only if $r_{x_1}<l_{x_2}$. We call $I$ an interval representation of $(X,P)$. An interval order $(X,P)$ is a semi-order if there is an interval representation $I$ assigning to each element $x\in X$ a closed unit-length interval $I(x)=[r_x-1,r_x]$ on the real number line so that for all $x_1,x_2\in X$ we have $x_1<x_2$ if and only if $r_{x_1}<r_{x_2}-1$. Moreover, a proper interval representation is an interval representation where no interval is contained in the the interior of any other. It can be shown that an interval order $(X,P)$ is a semi-order if and only if it has a proper interval representation.

\subsection{On-Line Width of Interval Orders}

The on-line width $\text{olwi}(w)$ of the class of interval orders of width at most $w$ is the largest integer $k$ for which there exists a strategy that forces any algorithm to use $k$ chains to partition an interval order of width $w$. This variant was solved by Kierstead and Trotter in the early $80$'s.

\begin{theorem}[Kierstead and Trotter \cite{kie-tro-81}]
    $\emph{olwi}(w)=3w-2$.
\end{theorem}

Baier, Bosek and Micek \cite{bai-bos-07} showed that if every element presented is required to be a maximal element, then the on-line width is $2w-1$.

\subsection{On-Line Width of Interval Orders with Representation}

The on-line width $\text{olwi}_R(w)$ of the class of interval orders of width at most $w$ with representation is the largest integer $k$ for which there exists a strategy that forces any algorithm to use $k$ chains to partition an interval order of width $w$ presented as intervals. This means that instead of presenting the elements of the interval order as points, the elements are presented as intervals. These intervals provide an interval representation for a unique poset $(X,P)$. This variant of the problem was solved by Chrobak and \'Slusarek.

\begin{theorem}[Chrobak and \'Slusarek \cite{chr-slu-88}]
    $\olwi_R(w)=3w-2$.
\end{theorem} 

\subsection{First-Fit}

First-Fit may be the first on-line algorithm that comes to mind for partitioning posets into a minimal number of chains. However, First-Fit does not perform well for the class of all posets. Kierstead \cite{kie-86} proved that First-Fit could be forced to use arbitrarily many chains to partition a poset of width $2$. Bosek, Krawczyk and Szczypka \cite{bos-kra-10} showed that First-Fit provides a polynomial bound for the class of posets which do not induce two long incomparable chains. Let $\text{FF}(w)$ be the largest integer for which First-Fit can be forced to use $w$ chains to partition an interval order of width $w$. Kierstead \cite{kie-88} showed that $\text{FF}(w)\leq 40w$. Kierstead and Qin \cite{kie-qin-95} improved this to $25.72w$, and Pemmaraju, Raman, and Varadarajan \cite{pem-ram-04} improved it further to $10w$. Narayanaswamy and Subhash Babu \cite{nar-bab-08} noticed that the technique in \cite{pem-ram-04} yields $8w$. On the other hand, Kierstead, Smith and Trotter \cite{kie-smi-16} showed that for every $\epsilon>0$, $\text{FF}(w)>(5-\epsilon)w$ for sufficiently large $w$. The efficiency of First-Fit on intervals orders remains unknown. Nevertheless, Chrobak and \'Slusarek \cite{chr-slu-88} proved that the efficiency of First-Fit on semi-orders is $2w-1$. Over $30$ years later, it remains unknown whether a more efficient on-line algorithm exists for partitioning semi-orders presented with representation. 

\subsection{On-Line Width of Semi-Orders} 

The on-line width $\text{olws}(w)$ of the class of semi-orders of width at most $w$ is the largest integer $k$ for which there exists a strategy that forces any algorithm to use $k$ chains to partition a semi-order of width $w$. Chrobak and \'Slusarek showed that First-Fit needs at most $2w-1$ chains. Bosek, Felsner, Kloch, Krawczyk, Matecki and Micek later showed that any on-line algorithm can be forced to use $2w-1$ chains.

\begin{theorem}[Chrobak and \'Slusarek \cite{chr-slu-88}, Bosek et al. \cite{survey}]
    $\olws(w)=2w-1$.
\end{theorem}

Felsner, Kloch, Matecki, and Micek \cite{fel-klo-13} showed that if every element presented is required to be a maximal element, then the on-line width is $\lfloor{\frac{1+\sqrt{5}}{2}w\rfloor}$.

\subsection{On-Line Width of Semi-Orders with Representation}

The on-line width $\text{olws}_R(w)$ of the class of semi-orders of width at most $w$ with representation is the largest integer $k$ for which there exists a strategy that forces any algorithm to use $k$ chains to partition a semi-order of width $w$ presented as unit-intervals. As before, the semi-order is presented in the form of intervals instead of points, however, the intervals in this variant must all have length $1$. 

This problem was first considered by Chrobak and \'Slusarek \cite{chr-slu-88}. They showed that First-Fit needs at most $2w-1$ chains to partition a semi-order with representation. They also showed that any greedy algorithm can be forced to use $2w-1$ chains. However, it remained unkown whether a more optimal algorithm exists. Epstein and Levy \cite{eps-lev-05} constructed a strategy which forces $3k$ chains on a semi-order of width $2k$ presented with representation for any positive integer $k$. This provides the following previously best known bounds.
\[
    \lfloor\frac{3}{2}w\rfloor\leq\olws_R(w)\leq 2w-1
\]
Note that there is another potential candidate for the representation of a semi-order, i.e., a proper interval representation. Presenting proper intervals instead of unit-intervals may make a difference in the problem. In particular, presenting proper intervals may allow more chains to be forced. Nevertheless, since every unit-interval representation is proper, our bound holds true for both choices of representation.
While progress continues to be made for the general problem, as well as other variants, no improvements have been made to these bounds for almost $20$ years. In this paper, we improve the lower bound slightly by presenting a strategy which forces $3k+2$ chains on a semi-order of width $2k+1$ presented with representation for any positive integer $k$.

\section{Proof of Theorem}
Since in this variant, we introduce the elements of the poset $(X,P)$ as unit-intervals, we may define each element by a real number $r_i$. More specifically, we define each element introduced by the right endpoint of the interval in the representation so that if we introduce the element $x_i$ as the unit-interval $[r_i-1,r_i]$, we simply define $x_i$ by $x_i = r_i$. Assume that $w=2k+1$ for some positive integer $k$. The strategy consists of $5$ stages.

\subsection{Stage 1}%-----------------------------------------------------------

We begin by introducing a stack of intervals $x_1,\ldots,x_k$ so that $x_i=0$ for $i\in\{1,\ldots,k\}$. Notice that the intervals in Stage $1$ form an antichain, and hence, must each be assigned a distinct chain. Let $A$ denote the set of chains $\{a_1,\ldots ,a_k\}$ used in Stage $1$.

\begin{figure}[h!]
    \centering
    \resizebox{4.75in}{!}{\begin{tikzpicture}

% a-stack
\draw (0,0) -- (4,0) node[above left]{$A$} -- (4,0.9) -- (0,0.9) -- cycle;

% a-brace
\draw (4.1,0) -- (4.2,0);
\draw (4.1,0.9) -- (4.2,0.9);
\draw (4.15,0) -- (4.15,0.45) node[right]{$k$} -- (4.15,0.9);

%left stack
\draw (-7,0) node[left]{\tiny$a_1$} -- (-3,0);
\draw (-7,0.2) node[left]{\tiny$a_2$} -- (-3,0.2);
\draw (-7,0.9) node[left]{\tiny$a_k$} -- (-3,0.9);

\node at (-5,0.55)[circle,fill,inner sep = 0.1pt]{};
\node at (-5,0.7)[circle,fill,inner sep = 0.1pt]{};
\node at (-5,0.4)[circle,fill,inner sep = 0.1pt]{};

\draw[implies-implies,double equal sign distance] (-2,0.5) -- (-1,0.5);

\end{tikzpicture}}
    \caption{Stage 1: forcing the first $k$ chains.}
\end{figure}
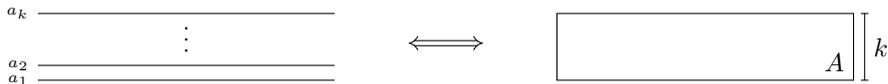

\subsection{Stage 2}%-----------------------------------------------------------

Initialize $l_2=1$ and $h_2=2$. In round $i$, we introduce the interval $x_i$ so that $x_i=(l_2+h_2)/2$. Suppose that the algorithm assigns the interval to chain $j$. If $j\in A$, then we update $h_2$ so that $h_2=x_i$. Otherwise, if $j\notin A$, then we update $l_2$ so that $l_2=x_i$. Let $B$ denote the set of new chains used in Stage $2$. If $|B|=k+1$, we move onto Stage $3$. Otherwise, if $|B|<k+1$, then we repeat Stage $2$ in round $i+1$.

Since $1<x_i<2$ for every interval $x_i$ introduced in Stage $2$, the intervals presented in Stage $2$ form an antichain of size at most $w$. Therefore, every interval is assigned to a different chain by the algorithm of which at most $k$ are in $A$. Hence, Stage $2$ ends forcing $k+1$ new chains. 

\begin{figure}[h!]
    \centering
    \resizebox{4.75in}{!}{\begin{tikzpicture}

\draw (0,0) -- (4,0) node[above left]{$A$} -- (4,0.9) -- (0,0.9) -- cycle;

\draw (5.5,0) -- (9.5,0) node[above left]{$B$} -- (10,1) -- (6,1) -- cycle;

\draw (7,1.1) -- (11,1.1) node[above left]{$\subseteq A$} -- (11.45,2) -- (7.45,2) -- cycle;

\draw[thin,dashed] (6,0) -- (6,2.2) node[above]{$l_2-1$};
    
\draw[thin,dashed] (7,0) -- (7,2.2) node[above]{$h_2-1$};

\draw (10.1,0) -- (10.2,0);

\draw (10.1,1) -- (10.2,1);

\draw (10.15,0) -- (10.15,0.5) node[right]{$k+1$} -- (10.15,1);
    
\end{tikzpicture}}
    \caption{Stage 2: forcing $k+1$ new chains.}
\end{figure}
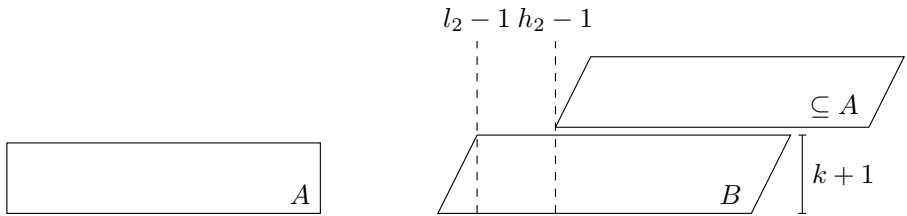

\subsection{Stage 3}%-----------------------------------------------------------

Initialize $l_3=l_2-3$ and $h_3=h_2-3$. In round $i$, we introduce a new interval $x_i$ so that $x_i=(l_3+h_3)/2$. Suppose that the algorithm assigns the interval to chain $j$. If $j\in B$, then we update $h_3$ so that $h_3=x_i$ and move onto Stage $4$. Otherwise, if $j\notin B$, then we update $l_3$ so that $l_3=x_i$ and we repeat Stage $3$ in round $i+1$ 

Since $-2<x_i<-1$ for every interval $x_i$ introduced in Stage $3$, the intervals presented in Stage $3$ form an antichain of size at most $w$. Therefore, every interval is assigned to a different chain by the algorithm of which at most $k$ are in $A$. If $k+1$ intervals are assigned entirely new chains, then the proof is complete. Hence, we may assume that Stage $3$ ends with the algorithm assigning an interval $x_B$ to a chain $b\in B$. Note that in this case, $x_B=h_3$

\begin{figure}[h!]
    \centering
    \resizebox{4.75in}{!}{\begin{tikzpicture}

\draw (0,0) -- (4,0) node[above left]{$A$} -- (4,0.9) -- (0,0.9) -- cycle;

%\draw (5.5,0) -- (9.5,0) node[above left]{$B$} -- (10,1) -- (6,1) -- cycle;

%\draw (7,1.1) -- (11,1.1) node[above left]{$\subset A$} -- (11.5,2.1) -- (7.5,2.1) -- cycle;

\draw (-1.9,0) node[above left]{$\notin B$} -- (-1.5,1.9) -- (-5.5,1.9) -- (-5.9,0) -- cycle;
    
\draw (-5.1,2) -- (-1.1,2) node[right]{$b$};

%\draw (-1.5,0.1) -- (-1.9,2.1) -- (-5.9,2.1) -- (-5.5,0.1) -- cycle;
    
%\draw (-5.1,0) -- (-1.1,0);
    
%\draw (-1.3,1.1) -- (-1.4,2.1) -- (2.6,2.1) -- (2.7,1.1) -- cycle;
    
%\draw (-1.2,1) -- (2.8,1);
    
%\draw (2.8,1.1) -- (6.8,1.1) -- (6.8,2.1) -- (2.8,2.1) -- cycle;
    
%\draw (-6,-1) -- (12,-1);
    
%\draw[thin,dashed] (6,-1) -- (6,3.2);
%\draw[thin,dashed] (7,-1) -- (7,3.2);
%\draw[thin,dashed] (2,-1) -- (2,3.2);
    
%\draw[thin,dashed] (3,-1) -- (3,3.2);
%\draw[thin,dashed] (-2,-1) -- (-2,3.2);
%\draw[thin,dashed] (-1,-1) -- (-1,3.2);
    
\draw[thin,dashed] (-1.5,0) -- (-1.5,2.3) node[above]{$l_3$};
%\draw[very thin,dashed] (2.5,-1) -- (2.5,3.2);
%\draw[very thin,dashed] (6.5,-1) -- (6.5,3.2);

\draw[thin,dashed] (-1.1,0) -- (-1.1,2.3) node[above]{$h_3$};
%\draw[very thin,dashed] (2.9,-1) -- (2.9,3.2);
%\draw[very thin,dashed] (6.9,-1) -- (6.9,3.2);
    
\end{tikzpicture}}
    \caption{Stage $3$: forcing a chain $b\in B$ on $x_B$.}
\end{figure}
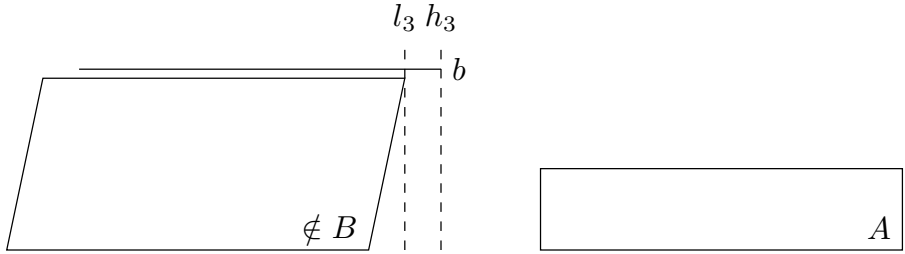

\subsection{Stage 4}%-----------------------------------------------------------

Initialize $l_4=l_3+1$ and $h_4=h_3+1$. In round $i$, we introduce a new interval $x_i$ so that $x_i=(l_4+h_4)/2$.  Suppose that the algorithm assigns the interval to chain $j$. Since $-1<x_i<x_b+1<0$, $j\notin A$ and $j\neq b$. We update $l_4$ so that $l_4=x_i$. If $j \notin B$, then we move onto Stage $5$. Otherwise, if $j\in B$, we repeat Stage $4$ in round $i+1$.

The intervals introduced in Stage $4$ form an antichain of size at most $k+1$. Therefore, every interval is assigned to a distinct chain of which no chain is in $A\cup\{b\}$ and at most $k$ chains are in $B\setminus \{b\}$. Hence, Stage $4$ ends with the algorithm assigning an interval $x_C$ to an entirely new chain $c$.

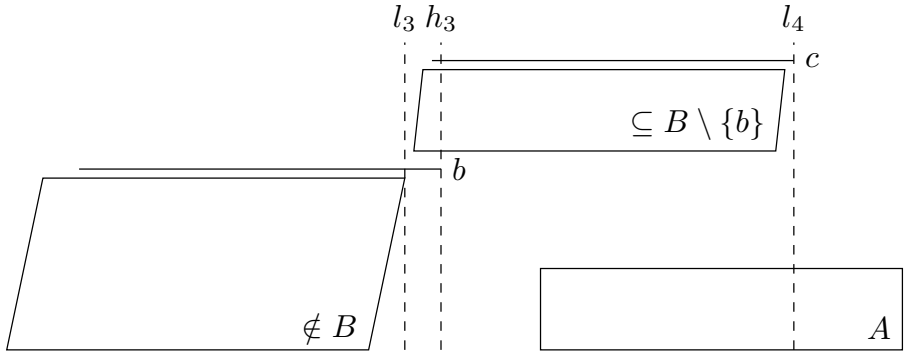
\begin{figure}[h!]
    \centering
    \resizebox{4.75in}{!}{\begin{tikzpicture}

\draw (0,0) -- (4,0) node[above left]{$A$} -- (4,0.9) -- (0,0.9) -- cycle;

%\draw (5.5,0) -- (9.5,0) node[above left]{$B$} -- (10,1) -- (6,1) -- cycle;
%\draw (7,1.1) -- (11,1.1) node[above left]{$\subset A$} -- (11.5,2.1) -- (7.5,2.1) -- cycle;

\draw (-1.9,0) node[above left]{$\notin B$} -- (-1.5,1.9) -- (-5.5,1.9) -- (-5.9,0) -- cycle;
    
\draw (-5.1,2) -- (-1.1,2) node[right]{$b$};

%\draw (-1.5,0.1) -- (-1.9,2.1) -- (-5.9,2.1) -- (-5.5,0.1) -- cycle;
%\draw (-5.1,0) -- (-1.1,0);

\draw (-1.4,2.2) -- (-1.3,3.1) -- (2.7,3.1) -- (2.6,2.2) node[above left]{$\subseteq B\setminus\{b\}$}-- cycle;
    
\draw (-1.2,3.2) -- (2.8,3.2) node[right]{$c$};
    
%\draw (-1.3,1.1) -- (-1.4,2.1) -- (2.6,2.1) -- (2.7,1.1) -- cycle;    
%\draw (-1.2,1) -- (2.8,1);    
%\draw (2.8,1.1) -- (6.8,1.1) -- (6.8,2.1) -- (2.8,2.1) -- cycle;    
%\draw (-6,-1) -- (12,-1);
    
%\draw[thin,dashed] (6,-1) -- (6,3.2);
%\draw[thin,dashed] (7,-1) -- (7,3.2);
%\draw[thin,dashed] (2,-1) -- (2,3.2);
    
%\draw[thin,dashed] (3,-1) -- (3,3.2);
%\draw[thin,dashed] (-2,-1) -- (-2,3.2);
%\draw[thin,dashed] (-1,-1) -- (-1,3.2);
    
\draw[thin,dashed] (-1.5,0) -- (-1.5,3.4) node[above]{$l_3$};
%\draw[very thin,dashed] (2.5,-1) -- (2.5,3.2);
%\draw[very thin,dashed] (6.5,-1) -- (6.5,3.2);

\draw[thin,dashed] (-1.1,0) -- (-1.1,3.4) node[above]{$h_3$};
\draw[thin,dashed] (2.8,0) -- (2.8,3.4) node[above]{$l_4$};
%\draw[very thin,dashed] (6.9,-1) -- (6.9,3.2);
    
\end{tikzpicture}}
    \caption{Stage $4$: forcing a new chain $c$ on $x_C$.}
\end{figure}

\subsection{Stage 5}%-----------------------------------------------------------

Finally, for each $i\in\{c+1,\dots,c+k\}$, we introduce an interval $x_i$ so that $x_i=x_C+1$. The intervals introduced in Stage $5$ form an antichain of size $k$ of which each interval cannot be assigned to any chain in $A\cup B\cup \{c\}$. All that is left to show is that we have not exceeded the width $w$. Let $x$ be any interval introduced in Stage $5$. It is trivial to check that the only interval from Stages $3$ and Stage $4$ which is incomparable to $x_i$ is $x_C$. Moreover, solving for the following: 
\[
l_2 -3 < x_B < h_2 -3\]
\[
l_2 -2 < x_C < x_B +1\]
\[
x_i = x_C +1\]
we get that $l_2-1<x_i<h_2-1$ which implies that the only intervals from Stage $2$ that are incomparable to $x_i$ are exactly the $k+1$ intervals which were assigned to chains from $B$. Let $D$ denote the set of new chains forced in Stage $5$. Thus, the total number of chains forced on this poset of width $w$ is
\[
|A|+|B|+|\{c\}|+|D|=k+(k+1)+1+k= 3k+2.
\]

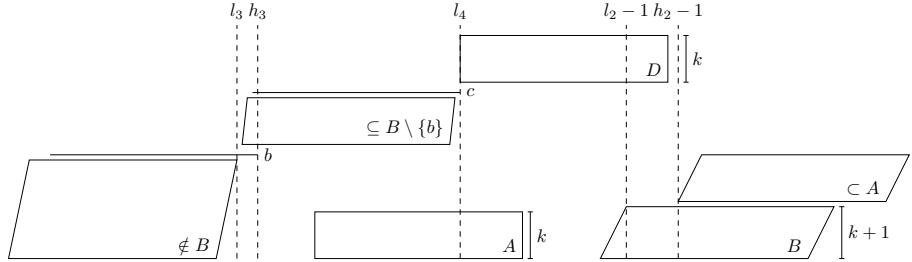
\begin{figure}[h!]
    \centering
    \resizebox{4.75in}{!}{\begin{tikzpicture}

\draw (0,0) -- (4,0) node[above left]{$A$} -- (4,0.9) -- (0,0.9) -- cycle;

\draw (5.5,0) -- (9.5,0) node[above left]{$B$} -- (10,1) -- (6,1) -- cycle;

\draw (7,1.1) -- (11,1.1) node[above left]{$\subset A$} -- (11.45,2) -- (7.45,2) -- cycle;

\draw (-1.9,0) node[above left]{$\notin B$} -- (-1.5,1.9) -- (-5.5,1.9) -- (-5.9,0) -- cycle;
    
\draw (-5.1,2) -- (-1.1,2) node[right]{$b$};

%\draw (-1.5,0.1) -- (-1.9,2.1) -- (-5.9,2.1) -- (-5.5,0.1) -- cycle;
    
%\draw (-5.1,0) -- (-1.1,0);

\draw (-1.4,2.2) -- (-1.3,3.1) -- (2.7,3.1) -- (2.6,2.2) node[above left]{$\subseteq B\setminus\{b\}$}-- cycle;
    
\draw (-1.2,3.2) -- (2.8,3.2) node[right]{$c$};
    
%\draw (-1.3,1.1) -- (-1.4,2.1) -- (2.6,2.1) -- (2.7,1.1) -- cycle;
    
%\draw (-1.2,1) -- (2.8,1);
    
\draw (2.8,3.4) -- (6.8,3.4) node[above left]{$D$} -- (6.8,4.3) -- (2.8,4.3) -- cycle;
    
%\draw (-6,-1) -- (12,-1);
    
%\draw[thin,dashed] (6,-1) -- (6,3.2);
%\draw[thin,dashed] (7,-1) -- (7,3.2);
%\draw[thin,dashed] (2,-1) -- (2,3.2);
    
%\draw[thin,dashed] (3,-1) -- (3,3.2);
%\draw[thin,dashed] (-2,-1) -- (-2,3.2);
%\draw[thin,dashed] (-1,-1) -- (-1,3.2);
    
\draw[thin,dashed] (-1.5,0) -- (-1.5,4.5) node[above]{$l_3$};
%\draw[very thin,dashed] (2.5,-1) -- (2.5,3.2);
%\draw[very thin,dashed] (6.5,-1) -- (6.5,3.2);

\draw[thin,dashed] (-1.1,0) -- (-1.1,4.5) node[above]{$h_3$};
\draw[thin,dashed] (2.8,0) -- (2.8,4.5) node[above]{$l_4$};
%\draw[very thin,dashed] (6.9,-1) -- (6.9,3.2);
    
\draw[thin,dashed] (6,0) -- (6,4.5) node[above]{$l_2-1$};
    
\draw[thin,dashed] (7,0) -- (7,4.5) node[above]{$h_2-1$};

\draw (4.1,0) -- (4.2,0);

\draw (4.1,0.9) -- (4.2,0.9);

\draw (4.15,0) -- (4.15,0.45) node[right]{$k$} -- (4.15,0.9);

\draw (10.1,0) -- (10.2,0);

\draw (10.1,1) -- (10.2,1);

\draw (10.15,0) -- (10.15,0.5) node[right]{$k+1$} -- (10.15,1);

\draw (7.1,3.4) -- (7.2,3.4);

\draw (7.1,4.3) -- (7.2,4.3);

\draw (7.15,3.4) -- (7.15,3.85) node[right]{$k$} -- (7.15,4.3);

\end{tikzpicture}}
    \caption{Stage 5: forcing the last $k$ chains.}
\end{figure}
This concludes the proof. 

\subsection{Remarks}

We proved that if $w=2k+1$, then our strategy will force any on-line algorithm to use $3k+2$ chains for any positive integer $k$. Moreover, we know that any greedy algorithm uses at most $2w-1$ chains. Thus, we get the answer to the previously open problem of finding the on-line width of the class of semi-orders of width $3$ with representation.

\begin{corr}
$\olws_R(3)=5$.
\end{corr}

\bibliographystyle{acm}
\bibliography{semi}

\end{document}